\documentclass[11pt]{amsart}

\usepackage{amsfonts,amsmath,amssymb,amsthm}
\usepackage[all]{xypic}

\newcommand{\Cc}{\mathbb{C}} 

\newcommand{\Pp}{\mathbb{P}}
\newcommand{\Xx}{\mathbb{X}}

\newcommand{\defi}[1]{\emph{#1}}

\newcommand{\Baff}{{\mathcal{B}_\mathit{\!aff}}}
\newcommand{\Binf}{{\mathcal{B}_\mathit{\!inf}}}
\newcommand{\B}{{\mathcal{B}}}
\newcommand{\gen}{\mathit{gen}}
\newcommand{\id}{\mathop{\mathrm{id}}\nolimits}
\newcommand{\cst}{\mathrm{cst}}
\newcommand{\PP}{{\mathcal{P}}}
\newcommand{\Sing}{\mathop{\mathrm{Sing}}\nolimits}

\renewcommand{\epsilon}{\varepsilon}
\renewcommand{\le}{\leqslant}
\renewcommand{\ge}{\geqslant}
\renewcommand {\leq}{\leqslant}

\newtheorem{theorem}{Theorem}

\newtheorem{lemma}[theorem]{Lemma}
\newtheorem{proposition}[theorem]{Proposition}

\theoremstyle{definition}
\newtheorem{remark}[theorem]{Remark}
\newtheorem{definition}[theorem]{Definition}
\newtheorem{example}[theorem]{Example}

\numberwithin{equation}{section}

\begin{document}
\title{Topological equivalence of complex polynomials}

\author{Arnaud Bodin}
\email{Arnaud.Bodin@math.univ-lille1.fr}
\author{Mihai Tib\u ar}
\email{Mihai.Tibar@math.univ-lille1.fr}
\address{Laboratoire Paul Painlev\'e, Math\' ematiques,
Universit\'e  Lille I, 
 59655 Villeneuve d'Ascq, France.}


\begin{abstract}
  The following numerical control over the topological equivalence is
  proved: two complex polynomials in $n\not= 3$ variables and with isolated
  singularities are
  topologically equivalent if one deforms into the other by a
  continuous family of polynomial functions $f_s \colon \Cc^n \to \Cc$
  with isolated singularities such that the degree, the number of
  vanishing cycles and the number of atypical values are constant in
  the family.
\end{abstract}

\maketitle


\section{Introduction}

Two polynomial functions $f,g : \Cc^n \to \Cc$ are said to be
\defi{topologically equivalent} if there exist homeomorphisms $\Phi :
\Cc^n \to \Cc^n$ and $\Psi:\Cc \to \Cc$ such that $\Psi \circ f =
g\circ \Phi$. A challenging natural question is: {\em under what
  conditions this topological equivalence is controlled by numerical
  invariants?}

We shall assume that our polynomials have isolated critical points,
and therefore finitely many.  It appears that the topology of a
polynomial function depends not only on critical points but also,
since the function is non-proper, on the behaviour of its fibres in
the neighbourhood of infinity.  It is well-known (and goes back to
Thom \cite{Th}) that a polynomial function has a finite set of
atypical values $\B \subset \Cc$, i.e. values at which the function
fails to be a locally trivial fibration.  Even if the critical points
of the polynomial are isolated, the homology of fibres may be very
complicated, behaving as if highly non-isolated singularities occur at
infinity.  One has studied such kind of singularities at infinity in
case they are in a certain sense isolated, e.g. \cite{Br}, \cite{ST}, \cite{LT},
  \cite{ST-gendefo}.  In this case the reduced homology of the general fibre
$G$ is concentrated in dimension $n-1$ and certain numbers may
be attached to singular points at infinity.

Coming back to topological equivalence: if our $f$ and $g$ are
topologically equivalent then clearly their corresponding fibres
(general or atypical) are homeomorphic. In particular the Euler
characteristics of the general fibres of $f$ and $g$ and the numbers
of atypical values of $f$ and $g$ coincide respectively.  We prove the
following numerical criterion for topological equivalence (see \S
\ref{sec:def} for an example):
\begin{theorem}\label{th:mucst}  
  Let $(f_s)_{s\in[0,1]}$ be a continuous family of complex
  polynomials with isolated singularities in the affine space and at
  infinity, in $n\not= 3$ variables.  If the numbers $\chi(G(s))$,
  $\#\B(s)$ and $\deg f_s$ are independent of $s\in[0,1]$, then the
  polynomials $f_0$ and $f_1$ are topologically equivalent.
\end{theorem}

In case of a smooth family of germs of holomorphic functions with
isolated singularity $g_s : (\Cc^n,0) \to (\Cc, 0)$, a famous result
by L\^{e}~D.T. and C.P.~Ramanujam \cite{LR} says that the constancy of
the {\em local Milnor number} (equivalently, of the Euler
characteristic of the general fibre in the local Milnor fibration
\cite{Mi}) implies that the hypersurface germs $g_0^{-1}(0)$ and
$g_1^{-1}(0)$ have the same topological type whenever $n\not= 3$.
J.G.~Timourian \cite{Tim} and H.~King \cite{Ki} showed moreover the
topological triviality of the family of function germs over a small
enough interval.  The techniques which are by now available for
proving the L\^e-Ramanujam-Timourian-King theorem do not work beyond
the case of isolated singularities. In other words, the topological
equivalence problem is still unsolved for local non-isolated
singularities.

The global setting poses new problems since one has to deal in the
same time with several singular points and atypical values.
Singularities at infinity introduce a new and essential difficulty
since they are of a different type than the critical points of
holomorphic germs.  Some evidence for the crucial importance of
singularities at infinity, even when assumed isolated, in
understanding the behaviour of polynomials is the famous unsolved
Jacobian Conjecture. One of the equivalent formulations of this
conjecture, in two variables, is the following, see \cite{LW}, \cite{ST}: if
$f : \Cc^2 \to \Cc$ has no critical points but has isolated
singularities at infinity then, for any polynomial $h : \Cc^2 \to
\Cc$, the critical locus $Z(\mathrm{Jac}(f,h))$ is not empty.

\smallskip

Our approach consists in three main steps, which we briefly describe in
the following.

\noindent 
{\em Step 1}. We show how the assumed invariance of numbers implies
the rigidity of the singularities, especially of those at infinity. We
use the following specific definition of {\em isolated singularities
  at infinity}, also employed in other papers (see e.g.  Libgober
\cite{Li} and Siersma-Tib\u ar \cite{ST-gendefo}): the projective
closure of any fibre of the polynomial and its slice by the hyperplane
at infinity have isolated singularities.  We explain in \S
\ref{sec:fam} how this condition enters in the proof of the key
results on the semi-continuity of certain numbers, which in turn
implies the rigidity of singularities.  The class of polynomials with
isolated singularities at infinity is large enough to include all
polynomials of two variables with reduced fibres.

\smallskip
\noindent 
{\em Step 2}. Assuming the rigidity of some singularity, we prove a
local topological triviality result. In case of a singular point in
$\Cc^n$ this is Timourian's theorem, but we have to consider the new
case of a singularity at infinity.  To handle such a problem we first
compactify all the polynomials of the family in the ``same way'', i.e.
we consider the total space $\Xx$ of a family depending polynomially
on the parameter: here we need the constancy of the degree.  One
cannot conclude by Thom-Mather's first Isotopy Lemma since the natural
stratification given by the trajectory of the singular point is not
Whitney in general. Unlike the local case, in our setting the
underlying space $\Xx$ turns out to be singular (essentially since the
compactification has singularities at infinity). Our strategy is to
revisit and modify the explicit local trivialisation given by
Timourian's proof by taking into account the Whitney stratification of
$\Xx$ (\S \ref{sec:locinf}).  The use of Timourian's proof is also
responsible for the excepted case $n= 3$, due to an argument by
L\^e-Ramanujam which relies on the h-cobordism theorem, cf \cite{LR}.

\smallskip
\noindent 
{\em Step 3}. Finally, we show how to patch together all the pieces
(i.e.  some open subsets of $\Cc^n$) where we have topological
triviality, in order to obtain the global topological equivalence. The
first named author used patching in \cite{Bo} to prove topological
equivalence in case there are no singularities at infinity and in case
$n=2$ with additional hypotheses and relying on results by L. Fourrier
\cite{Fo} which involve resolution of singularities. In our setting we
have to deal with pieces coming from singularities at infinity and
their patching is more delicate (see \S \ref{sec:proof}).
 
Let us remark that our theorem only requires the continuity of the
family instead of the smoothness in \cite{LR}, \cite{Tim}.  The reduction
from a continuous family to a family depending polynomially on the
parameter is made possible by a constructibility argument developed in
\S \ref{s:construct}.  The constructibility also implies the
finiteness of topological types of polynomials when fixing numerical
invariants, see Remark \ref{r:finite}.  It is worth to point out that
the finiteness does not hold for the equivalence up to
diffeomorphisms, as already remarked by T. Fukuda \cite{Fu}.  For
example, the family $f_s(x,y) = xy(x-y)(x-sy)$ provides infinitely
many classes for this equivalence, because of the variation of the
cross-ratio of the $4$ lines.

\section{Definitions and notations}
\label{sec:def}

We consider a one-parameter family of polynomials $f_s(x) = P(x,s)$,
where $P:\Cc^n \times [0,1] \to \Cc$ is polynomial in $s$ and such
that $\deg f_s =d$, for all $s\in [0,1]$.

We assume that the \defi{affine singularities of $f_s$ are isolated}:
$\dim \Sing f_s \le 0$ for all $s$, where $\Sing f_s = \{ x \in \Cc^n
\mid \mathrm{grad}\, f_s (x) = 0\}$. The set of \defi{affine critical values} of
$f_s$ is a finite set and we denote it by $\Baff(s) = \{ t \in \Cc
\mid \mu_t(s) > 0\}$, where $\mu_t(s)$ is the sum of the local Milnor
numbers at the singular points of the fibre $f_s^{-1}(t)$; remark that
we also have $\Baff(s) = f_s(\Sing f_s)$.  The \defi{total Milnor
  number} is $\mu(s) = \sum_{t\in \Baff} \mu_t(s)$.  We also assume
that, for all $s$, $f_s$ has isolated singularities at infinity in the
following sense.
\begin{definition}\label{d:isol}
  We say that a polynomial $f_s$ has \defi{isolated singularities at
    infinity} if $\dim \Sing W(s) \le 0$, where
\[
W(s) = \bigg\lbrace [x] \in \Pp^{n-1} \mid \frac{\partial
  P_d}{\partial x_1}=\cdots =\frac{\partial P_d}{\partial x_n} =0
\bigg\rbrace
\]
is an algebraic subset of the hyperplane at infinity $H^\infty$ of
$\Pp^n$, which we identify to $\Pp^{n-1}$. Here $P_d$ denotes the
homogeneous part of degree $d$ in variables $x=(x_1,\ldots,x_n)$ of
the polynomial $P(x,s)$.  The condition $\Sing W(s) \le 0$ is
equivalent to the following: for all $t\in \Cc$, the singularities of
$\overline{f_s^{-1}(t)}$ and of $\overline{f_s^{-1}(t)}\cap H^\infty$
are at most isolated.
\end{definition}
In \cite{ST-gendefo} one calls ``$\mathcal F$ISI  deformation of $f_0$'' a
family $P$ such that $f_s$ has isolated singularities at infinity and
in the affine, for all $s$.  The class of polynomials with isolated
singularities at infinity is large enough to include all polynomial
functions in two variables with reduced fibres.  It is a (strict)
subclass of polynomials having isolated singularities at infinity in
the sense used by Broughton \cite{Br} or in the more general sense of
\cite{ST}.

We shall see in the following how one can precisely detect the
singularities at infinity.  We attach to the family $P$ the following
hypersurface:
 \[ \Xx = \big\lbrace{ ([x:x_0], t, s) \in \Pp^n \times \Cc \times [0,1] \mid \tilde 
   P(x,x_0,s) - tx_0^d = 0 \big\rbrace},\] where $\tilde P$ denotes
 the homogenisation of $P$ by the new variable $x_0$, considering $s$
 as parameter.  Let $\tau : \Xx \to \Cc$ be the projection to the
 $t$-coordinate. This extends the map $P$ to a proper one in the sense
 that $\Cc^n\times [0,1]$ is embedded into $\Xx$ and that the
 restriction of $\tau$ to $\Cc^n\times [0,1]$ is equal to $P$.  Let
 $\sigma : \Xx \to [0,1]$ denote the projection to the $s$-coordinate.
 We shall use the notations $\Xx(s) = \sigma^{-1}(s) \cap \Xx$.  Let
 $\Xx_{t}(s) = \Xx(s) \cap \tau^{-1}(t)$ be the projective closure in
 $\Pp^n$ of the affine hypersurface $f_s^{-1}(t)$.  Note that $\Xx(s)$
 is singular, in general with 1-dimensional singular locus, since
 $\Sing \Xx(s) = \Sigma(s) \times \Cc$, where:
 \[
 \Sigma(s) = \bigg\lbrace [x] \in \Pp^{n-1} \mid \frac{\partial
   P_d}{\partial x_1}(x,s)=\cdots =\frac{\partial P_d}{\partial
   x_n}(x,s)= P_{d-1}(x,s) =0 \bigg\rbrace
\]
and we have $\Sigma(s) \subset W(s)$, which implies that $\Sigma(s)$
is finite.

Let us fix some $s \in [0,1]$ and some $p\in\Sigma(s)$. For $t\in \Cc$,
let $\mu_p(\Xx_{t}(s))$ denote the local Milnor number of the
projective hypersurface $\Xx_{t}(s) \subset \Pp^n$ at the point
$[p:0]$.  By \cite{Br} the number $\mu_p(\Xx_{t}(s))$ is constant for
generic $t$, and we denote this value by $\mu_{p,\gen}(s)$.  We have
$\mu_p(\Xx_{t}(s)) > \mu_{p,\gen}(s)$ for a finite number of values of
$t$.  The \defi{Milnor-L\^{e} number} at the point $([p:0], t)\in
\Xx(s)$ is defined as the jump $\lambda_{p,t}(s) := \mu_p(\Xx_{t}(s))
- \mu_{p,\gen}(s)$.  We say that the point $([p:0], t)$ is a
\defi{singularity at infinity} of $f_s$ if $\lambda_{p,t}(s)> 0$.  Let
$\lambda_t(s) = \sum_{p\in \Sigma(s)}\lambda_{p,t}(s)$.  The set of
\defi{critical values at infinity} of the polynomial $f_s$ is defined
as:
$$\Binf(s) = \{ t\in \Cc \mid \lambda_t(s) > 0\}.$$
Finally, the
\defi{Milnor-L\^e number at infinity} of $f_s$ is defined as:
$$\lambda(s) = \sum_{t\in\Binf(s)}\lambda_t(s).$$

For such a polynomial, the \defi{set of atypical values}, or the
\defi{bifurcation set}, is:
$$\B(s) = \Baff(s) \cup \Binf(s).$$
It is known that $f_s :
f_s^{-1}(\Cc \setminus \B(s)) \to \Cc \setminus \B(s)$ is a locally
trivial fibration \cite{Br}.  After \cite{ST}, for $t\in \Cc$ the
fibre $f_s^{-1}(t)$ is homotopic to a wedge of spheres of real
dimension $n-1$ and the number of these spheres is
$\mu(s)+\lambda(s)-\mu_t(s)-\lambda_t(s)$.  In particular, for the
Euler characteristic of the general fibre $G(s)$ of $f_s$ one has:
\[ \chi(G(s)) = 1 + (-1)^{n-1}(\mu(s)+\lambda(s)).\]
\begin{example}
  Let $f_s(x,y,z,w) = x^2y^2+ z^2+ w^2 + xy+(1+s)x^2+x$.  For $s \in
  \Cc \setminus \{-2,-1\}$ we have $\B(s) = \big\lbrace 0, -\frac14,
  -\frac14\frac{s+2}{s+1}\big \rbrace$, $\mu(s) = 2$ and $\lambda(s)=
  1$.  It follows that $\chi(G(s)) =1- \mu(s)-\lambda(s) = -2$ and
  that $\# \B(s) = 3$. For the two excepted polynomials $f_{-1}$ and
  $f_{-2}$ we have $\#\B = 2$. Then, by Theorem \ref{th:mucst}, $f_0$
  is topologically equivalent to $f_s$ if and only if $s \in \Cc
  \setminus \{-2,-1\}$.
\end{example}
Let $\Xx^\infty(s)$ denote the part ``at infinity'' $\Xx(s)\cap \{ x_0
=0\}$. We shall use in \S\ref{sec:locinf} the Whitney stratification
of the space $\Xx(s)$ with the following strata (see \cite{ST},
 \cite{ST-gendefo}): $\Xx(s)\setminus\Xx^\infty(s)$, $\Xx^\infty(s)
\setminus\Sing \Xx(s)$, the complement in $\Sing \Xx(s)$ of the
singularities at infinity and the finite set of singular points at
infinity.  We also recall that the restriction $\tau : \Xx(s) \to \Cc$
is transversal to all the strata of $\Xx(s)$ except at the singular
points at infinity.

\section{Rigidity of singularities in families of polynomials}
\label{sec:fam}

Let $(f_s)_{s\in [0,1]}$ be a family of complex polynomials with
constant degree $d$, such that the coefficients of $f_s$ are
polynomial functions of $s \in [0,1]$. We also suppose that for all $s
\in [0,1]$, $f_s$ has isolated singularities in the affine space and
at infinity (in the sense of Definition \ref{d:isol}). Under these
conditions, we may prove the following rigidity result:

\begin{proposition}\label{p:rigid}
  If the pair of numbers $(\mu(s)+\lambda(s), \#\B(s))$ is independent
  of $s$ in some neighbourhood of $0$, then the $5$-uple $(\mu(s),
  \#\Baff(s), \lambda(s), \#\Binf(s),\#\B(s))$ is independent of $s$
  too.  Moreover there is no collision of points $p(s) \in \Sigma(s)$
  as $s\to 0$, and in particular $\# \Sigma(s)$ and $\mu_{p,\gen}(s)$
  are constant.
\end{proposition}

\begin{proof}
  {\bf Step 1.} We claim that the multivalued map $s\mapsto \B(s)$ is
  continuous. If not the case, then there is some value of $\B(s)$
  which disappears as $s\to 0$. To compensate this, since $\#\B(s)$ is
  constant, there must be a value which appears in $\B(0)$. By the
  local constancy of the total Minor number, affine singularities
  cannot appear from nothing, therefore the new critical value should
  be in $\Binf(0)$. More precisely, there is a singular point at
  infinity $(p,t)$ of $f_0$ (thus where the local $\lambda$ is
  positive) such that, for $s\not= 0$, there is no singular point of
  $f_s$, either in affine space or at infinity, which tends to $(p,t)$
  as $s\to 0$. But this situation has been shown to be impossible in
  \cite[Lemma 20]{Bo}.  Briefly, the argument goes as follows: Let
  $(p,c(0))$ be a singularity at infinity of $f_0$ and let
  $h_{s,t}:\Cc^n\rightarrow \Cc$ be the localisation at $p$ of the map
  $\tilde{P}(x,x_0,s)-tx_0^d$. Then from the local conservation of the
  total Milnor number of $h_{0,c(0)}$ and the dimension of the
  critical locus of the family $h_{s,t}$ one draws a contradiction.
  The claim is proved.
  
  Let us remark that our proof also implies that the finite set
  $\B(s)\subset \Cc$ is contained in some disk of radius independent
  of $s$ and that there is no collision of points of $\B(s)$ as
  $s\rightarrow 0$.

\noindent
{\bf Step 2.}  We prove that there is no collision of points $p(s) \in
\Sigma(s)$ as $s\to 0$ and that $\# \Sigma(s)$ and $\mu_{p,\gen}(s)$
are constant.  We pick up and fix a value $t\in \Cc$ such that
$t\not\in \B(s)$, for all $s$ near $0$.  Then we have a one parameter
family of general fibres $f_s^{-1}(t)$, where $s$ varies in a
neighbourhood of 0. The corresponding compactified hypersurfaces
$\Xx_{t}(s)$ have isolated singularities at their intersections with
the hyperplane at infinity $H^\infty$.

Let $\mu_p^\infty(s)$ denote the Milnor number of the hyperplane slice
$\Xx_{t}(s) \cap H^\infty$ at some $p\in W(s)$, and note that this
does not depend on $t$, for some fixed $s$.  We use the following
formula (see \cite[2.4]{ST-gendefo} for the proof and references):
\begin{equation}\label{form:equality}
 \mu(s)+\lambda(s) = (d-1)^n - \sum_{p\in\Sigma(s)}\mu_{p,\gen}(s)-\sum_{p\in
   W(s)}\mu_p^\infty(s).
\end{equation}
Since $\mu(s) + \lambda(s)$ is constant and since the local upper
semi-continuity of Milnor numbers, we have that both sums
$\sum_{p\in\Sigma(s)}\mu_{p,\gen}(s)$ and $\sum_{p\in W(s)}
\mu_p^\infty(s)$ are constant hence locally constant.  The
non-splitting principle (see \cite{La} or \cite{Le2}, \cite{AC}) applied to
our family of hypersurface multigerms tells that each
$\mu_{p,\gen}(s)$ has to be constant. This means that there cannot be
collision of points of $\Sigma(s)$.

\noindent
{\bf Step 3.} We claim that $\mu(s)$ is constant. If not the case,
then we may suppose that $\mu(0) < \mu(s)$, for $s$ close to 0, since
$\mu(s)$ is lower semi-continuous (see \cite{Br}).  Then by also using
Step 1, there exists $c(s) \in \Baff(s)$, such that: $c(s) \to c(0)
\in \Cc$ as $s \rightarrow 0$.  By Step 1, there is no other value
except $c(s)\in \B(s)$ which tends to $c(0)$. We therefore have a
family of hypersurfaces $\Xx_{c(s)}(s)$ with isolated singularities
$q_j(s) \in f_s^{-1}(c(s))$ that tend to the singularity at $(p,0)\in
\Sigma(0) \subset \Xx_{c(0)}(0)$.  By Step 2 and the (upper)
semi-continuity of the local Milnor numbers we have:

\begin{equation}\label{eq:mu1}
\mu_p(\Xx_{c(0)}(0)) \ge \mu_p(\Xx_{c(s)}(s)) + \sum_j \mu_{q_j(s)}(\Xx_{c(s)}(s)).
\end{equation}
By definition, $\mu_p(\Xx_{c(s)}(s)) = \lambda_p(s) + \mu_{p,\gen}(s)$
for any $s$, and by Step 2, $\mu_{p,\gen}(s)$ is independent of $s$.
It follows that:
\begin{equation}\label{eq:mu2}
 \lambda_{p, c(0)}(0) \ge \lambda_{p, c(s)}(s) + \sum_j\mu_{q_j(s)}(\Xx_{c(s)}(s)),
\end{equation}
which actually expresses the balance at any collision of singularities
at some point at infinity. This shows that in such collisions the
``total quantity of singularity'', i.e. the local $\mu + \lambda$, is
upper semi-continuous. On the other hand, the global $\mu + \lambda$
is assumed constant, by our hypothesis. This implies that the local
$\mu + \lambda$ is constant too. Therefore in (\ref{eq:mu2}) we must
have equality and consequently (\ref{eq:mu1}) is an equality too.

We may now conclude by applying the non-splitting principle, similarly
as in Step 2, to yield a contradiction.

\noindent
{\bf Step 4.} Since by Step 3 there is no loss of $\mu$, the
multi-valued function $s\mapsto \Baff(s)$ is continuous.  Steps 1 and
3 show that $s\mapsto \Binf(s)$ is continuous too.  Together with
$\#(\Baff(s) \cup \Binf(s)) = \cst$, this implies that
$\#\Baff(s)=\cst$ and $\#\Binf(s)=\cst$.
\end{proof}

\section{Constructibility via numerical invariants}
\label{s:construct}

Let $\PP_{\le d}$ be the vector space of all polynomials in $n$
complex variables of degree at most $d$.  We consider here the subset
$\PP_d(\mu+\lambda,\#\B)\subset \PP_{\le d}$ of polynomials of degree
$d$ with fixed $\mu+\lambda$ and fixed $\#\B$.

Recall that a \defi{locally closed set} is the intersection of a
Zariski closed set with a Zariski open set; a \defi{constructible set}
is a finite union of locally closed sets.

\begin{proposition}\label{l:construc}
  $\PP_d(\mu+\lambda,\#\B)$ is a constructible set.
\end{proposition}

\begin{proof}
  The set $\PP_d$ of polynomials of degree $d$ is a constructible set
  in the vector space $\PP_{\le d}$.  Let us first prove that
  ``isolated singularities at infinity'' yields a constructible set.
  A polynomial $f$ has isolated singularities at infinity if and only
  if $W:= W(f)$ has dimension $0$ or is void.  Let $S = \lbrace (x, f)
  \in \Pp^n \times \PP_d \mid f \in \PP_d,\ x\in W(f) \rbrace$ and let
  $\pi : S \to \PP_d$ be the projection on the second factor.  Since
  this is an algebraic map, by Chevalley's Theorem (e.g. \cite[\S
  14.3]{Eisenbud}) the set $\{ f \in \PP_d \mid \dim \pi^{-1}(f) \leq
  0\}$ is constructible and this is exactly the set of polynomials
  with isolated singularities at infinity.
  
  Next, we prove that fixing each integer $\mu$, $\#\Baff$, $\lambda$,
  $\#\Binf$, $\#\B$ yields a constructible set.  The main reason is
  the semi-continuity of the Milnor number (upper in the local case,
  lower in the affine case), see e.g. Broughton \cite[Prop. 2.3]{Br}.
  Broughton proved that the set of polynomials with a given $\mu <
  \infty$ is constructible.  As the inverse image of a constructible
  set by an algebraic map, the set of polynomials with Milnor number
  $\mu$ and bifurcation set such that $\#\Baff = k$ is a constructible
  set.
  
  Let $\PP_d(\mu,\#\Sigma)$ be the set of polynomials of degree $d$,
  with a given $\mu$, with isolated singularities at infinity and a
  given $\# \Sigma$. Notice that $\# \Sigma$ is finite because $\Sigma
  \subset W$ and is bounded for fixed $d$.  Since $\Sigma$ depends
  algebraically on $f$, we have that $\PP_d(\mu,\#\Sigma)$ is a
  constructible set.  Now the local Milnor number $\mu_p$ is an upper
  semi-continuous function, so fixing $\lambda_p$ as the difference of
  two local Milnor numbers (see \S \ref{sec:def}) provides a
  constructible set.  By doing this for all the critical points at
  infinity we get that fixing $\lambda = \sum_p \lambda_p$ yields a
  constructible condition.  The arguments for the conditions $\#\Binf$
  and $\#\B$ (which are numbers of points of two algebraic sets in
  $\Cc$) are similar to the one for $\#\Baff$.
  
  The just proved constructibility of $\PP_d(\mu, \#\Baff, \lambda,
  \#\Binf, \#\B)$ implies, by taking a finite union, the
  constructibility of $\PP_d(\mu+\lambda,\#\B)$.
\end{proof}

\begin{definition}
  We say that a finite set $\Omega(s)$ of points in $\Cc^k$, for some
  $k$, depending on a real parameter $s$, is an \defi{algebraic braid}
  if $\Omega= \cup_s \Omega(s)\times \{ s \}$ is a real algebraic
  sub-variety of $\Cc^k\times [0,1]$, the multi-valued function $s
  \mapsto \Omega(s)$ is continuous and  $\# \Omega(s) = \cst$.
\end{definition}

We may now reformulate and extend Proposition \ref{p:rigid} as
follows.
\begin{proposition}
\label{c:dep}
Let $(f_s)_{s\in[0,1]}$ be a family of complex polynomials with
isolated singularities in the affine space and at infinity, whose
coefficients are polynomial functions of $s$. Suppose that the numbers
$\mu(s)+\lambda(s)$, $\#\B(s)$ and $\deg f_s$ are independent of $s
\in[0,1]$.  Then:
\begin{enumerate}
\item \label{it:p1} $\Sigma(s)$, $\Baff(s)$, $\Binf(s)$ and $\B(s)$ are algebraic
  braids;
\item  \label{it:p2} for any continuous function $s\mapsto p(s) \in \Sigma(s)$ we
  have $\mu_{p(s),\gen}=\cst$;
\item  \label{it:p3} for any continuous function $s\mapsto c(s) \in \Binf(s)$ we have
  $\lambda_{p(s),c(s)}=\cst$;
\item  \label{it:p4} for any continuous function $s\mapsto c(s) \in \Baff(s)$ we have
  $\mu_{c(s)} = \cst$ and moreover, the local $\mu$'s of the fibre
  $f_s^{-1}(c(s))$ are constant.
\end{enumerate}
\end{proposition}

\begin{proof}
  (\ref{it:p1}) For $\Sigma(s)$, it follows from the algebraicity of the
  definition of $\Sigma$ and from Step 2 of Proposition \ref{p:rigid}.
  It is well-known that affine critical values of polynomials are
  algebraic functions of the coefficients.  Together with Proposition
  \ref{p:rigid}, this proves that $\Baff(s)$ is an algebraic braid.
  
  Similarly $\cup_s \Binf(s)\times \{ s \}$ is the image by a finite
  map of an algebraic set, and together with Step 4 of Proposition
  \ref{p:rigid}, this proves that $\Binf(s)$ is an algebraic braid.
   
  Next, (\ref{it:p2}) is Step 2 of Proposition \ref{p:rigid} and  (\ref{it:p3}) is a
  consequence of Step 3.  Lastly, observe that (\ref{it:p4}) is a well-known
  property of local isolated hypersurface singularities and follows
  from  (\ref{it:p1}) and the local non-splitting principle.
\end{proof}

\begin{remark} 
\label{r:finite}
Theorem \ref{th:mucst} has the following interpretation: to a
connected component of $\PP_d(\mu+\lambda, \#\B)$ one associates a
unique topological type. (It should be noticed that two different
connected components of $\PP_d(\mu+\lambda, \#\B)$ may have the same
topological type, see \cite{Bo2} for an example.)  
It follows that there is a finite number of topological types of complex
polynomials of fixed degree and with isolated singularities in the
affine space and at infinity. This may be related
to the finiteness of topological equivalence classes in $\PP_{\le d}$,
conjectured by Ren\'e Thom and proved by T.~Fukuda \cite{Fu}. 
\end{remark}
 
\section{Local triviality at infinity}
\label{sec:locinf}

The aim of this section is to prove a topological triviality statement
for a singularity at infinity. Our situation is new since it concerns
a family of couples space-function varying with the parameter $s$ and
where the space is singular. The proof actually relies on Timourian's
proof \cite{Tim} for germs of holomorphic functions on $\Cc^n$.
Therefore we shall point out where and how this proof needs to be
modified, since we have to plug-in a singular stratified space (i.e.
the germ of $\Xx(s)$ at a singularity at infinity) instead of the germ
$(\Cc^n,0)$ in Timourian's proof.

As before, let $(f_s)_{s\in[0,1]}$ be a family of complex polynomials
of degree $d$ and let $(p,c)$ be a singularity at infinity of $f_0$.
Let $g_s : \Xx(s) \to \Cc$ be the localisation at $(p(s),c(s))$ of the
map $\tau_{|\Xx(s)}$.  We denote by $B_\epsilon \subset \Cc^n\times
\Cc$ the closed $2n+2$-ball of radius $\epsilon$ centred at $(p,c)$,
such that $B_\epsilon\cap \Xx(0)$ is a Milnor ball for $g_0$.  We
choose $0 < \eta \ll \epsilon$ such that we get a Milnor tube $T_0 =
B_\epsilon \cap\Xx(0) \cap g_0^{-1}(D_\eta(c))$.  Then, for all $t \in
D_\eta(c)$, $g_0^{-1}(t)$ intersects transversally $S_\epsilon =
\partial B_\epsilon$.  We recall from \cite{ST} that $g_0 \colon
T_0\setminus g_0^{-1}(c) \to D_\eta(c)\setminus \{ c\}$ is a locally
trivial fibration whenever $\lambda_{p,c}(0) > 0$ and $g_0 \colon T_0
\to D_\eta(c)$ is a trivial fibration whenever $\lambda_{p,c}(0) = 0$.
 
According to Proposition \ref{c:dep}(\ref{it:p1}), by an analytic change of
coordinates, we may assume that $(p(s), c(s))=(p,c)$ for all $s\in
[0,u]$, for some small enough $u>0$.  We set $T_s = B_\epsilon\cap
\Xx(s) \cap g_s^{-1}(D_\eta(c))$ and notice that $B_\epsilon$ does not
necessarily define a Milnor ball for $g_s$ whenever $s\not= 0$.  For
some $u>0$, let $T = \bigcup_{s\in[0,u]} T_s \times \{s\}$, and let $G
: T \to \Cc \times [0,u]$ be defined by $G(z,s)=(g_s(z),s)$.

The homeomorphisms between the tubes that we consider here are all
{\em stratified}, sending strata to corresponding strata. The
stratification of some tube $T_s$ has by definition three strata: $\{
T_s \setminus (\{ p\} \times D_\eta(c)), \{ p\} \times
D_\eta(c)\setminus (p,c), (p,c)\}$.

\begin{theorem}
\label{th:muinf}
Let $f_s(x)= P(x,s)$ be a one-parameter polynomial family of
polynomial functions of constant degree, such that the numbers
$\mu(s)+\lambda(s)$ and $\# \B(s)$ are independent of $s$.  If $n\not=
3$, then there exists $u>0$ and a homeomorphism $\alpha$ such that the
following diagram commutes:
\[
\xymatrix{
  T \ar[r]^-\alpha \ar[d]_-{G}  & T_0 \times [0,u] \ar[d]^-{g_0 \times \id} \\
  D_\eta(c)  \times [0,u] \ar[r]_-\id   & D_\eta(c) \times [0,u], \\
}
\]
and such that $\alpha$ sends the strata of every $T_s$ to the
corresponding strata of $T_0$.
\end{theorem}
  
\begin{proof}
  Our point $(p,c)\in \Sigma(0)\times \Cc$ is such that
  $\lambda_{p,c}(0) > 0$. We cannot apply directly Timourian's result
  for the family $g_s$ because each function $g_s$ is defined on a
  {\em singular} space germ $(\Xx(s), (p,c))$, but we can adapt it to
  our situation. To do this we recall the main lines of this proof and
  show how to take into account the singularities via the
  stratification of $T_s$.

  Remark first that $(p,c)$ is the only singularity of $g_s$ in $T_s$,
  by the rigidity result Proposition \ref{c:dep}.  We use the notion
  of $\epsilon$-homeomorphism, meaning a homeomorphism which moves
  every point within a distance no more than $\epsilon >0$.
  
  Theorem \ref{th:muinf} will follow from Lemma \ref{lem:3tim} below,
  once we have proved that the assumptions of this lemma are fulfilled
  in our new setting. This is a simplified statement of Timourian's
  Lemma 3 in \cite{Tim}, its proof is purely topological and needs no
  change.

\begin{lemma} {\rm (\cite[Lemma 3]{Tim})}\label{lem:3tim} 
  Assume that:
\begin{enumerate}
\item \label{it:l1} The space of stratified homeomorphisms of $T_0$ into itself,
  preserving the fibres of $g_0$, is locally contractible.
\item  \label{it:l2} For any $\varepsilon >0$ there exists $u >0$ small enough such
  that for any $s, s'\in [0,u]$ there is a stratified
  $\varepsilon$-homeomorphism $h: T_s \to T_{s'}$ with $g_s = g_{s'}
  \circ h$.
\end{enumerate}
Then there exists a homeomorphism $\alpha$ as in Theorem
\ref{th:muinf}. 
\end{lemma}

The assumptions (\ref{it:l1}) and (\ref{it:l2}) correspond, respectively, to Lemma 1 and
Lemma 2 of Timourian's paper \cite{Tim}. The remaining proof is
therefore devoted to showing why the assumptions  (\ref{it:l1}) and (\ref{it:l2}) are true
in our setting.

Condition (\ref{it:l1}) can be proved as follows. It is well-known that analytic
sets have local conical structure \cite{BV}. Notice that the
stratification $\{ T_0 \setminus (\{ p\} \times D_\eta(c)), \{ p\}
\times D_\eta(c)\setminus (p,c), (p,c)\}$ of $T_0$ is a Whitney
stratification (but that this is not necessarily true for tubes $T_s$
with $s\not= 0$).  Timourian shows how to construct a vector field on
$T_0$ such that all integral curves end at the central point $(p,c)$.
Moreover, this vector field can be chosen such that to respect the
Whitney strata. This is the only new requirement that we need to plug
in. The rest of Timourian's argument remains unchanged once we have
got the vector field, and we give only its main lines in the
following.  This vector field is used to define a continuous family
$h_t$ of homeomorphisms such that $g_0 = g_0 \circ h_t$, which deforms
a homeomorphism $h_1 = h$ of $T_0$ which is the identity at the
boundary $\partial T_0$, to a homeomorphism $h_0$ which is the
identity within a neighbourhood of $(\partial B_\varepsilon \cap
T_0)\cup g_0^{-1}(c)\setminus (p,c)$. Next, by using the contracting
vector field, one constructs an isotopy of $h_0$ to the identity,
preserving the fibres of $g_0$.  To complete the proof, Timourian
shows how to get rid of the auxiliary condition ``to be the identity
at the boundary $\partial T_0$'' imposed to $h$, by using Siebenmann's
results \cite{Si}.

Condition (\ref{it:l2}) now. It will be sufficient to construct homeomorphisms
as in (\ref{it:l2}) from $T_0$ to $T_s$ for every $s\in [0,u]$, and take $u$
sufficiently small with respect to $\varepsilon$.  First remark that
for a sufficiently small $u$, the fibre $g_s^{-1}(t)$ intersects
transversally the sphere $S_\epsilon = \partial B_\epsilon$, for all
$s \in [0,u]$, and for all $t \in D_\eta(c)$. Consequently one may
define a homeomorphism $h' : \partial B_\epsilon \cap T_0 \to \partial
B_\epsilon \cap T_{s}$.  The problem is to extend it to an
homeomorphism from $T_0$ to $T_{s}$.

Take a Milnor ball $B_{\epsilon'}\subset B_\epsilon$ for $g_s$ at
$(p,c)$.  It appears that $(B_\epsilon\setminus \mathring
B_{\epsilon'})\cap g_s^{-1}(c)$ is diffeomorphic to $(\partial
B_\epsilon\cap g_s^{-1}(c)) \times [0,1]$.  This would be a
consequence of the h-cobordism theorem (the condition $n\not=3$ is
needed here) provided that it can be applied. The argument is given by
L\^e-Ramanujam's in \cite{LR} and we show how this adapts to our
setting.  One first notices that L\^e-Ramanujam's argument works as
soon as one has the following conditions: for $b\in D_\eta(c)$ and
$b\not= c$, the fibres $B_\epsilon\cap g_0^{-1}(b)$ and
$B_{\epsilon'}\cap g_s^{-1}(b)$ are singular only at $(p,b)$, they are
homotopy equivalent to a bouquet of spheres $S^{n-1}$ and the number
of spheres is the same. In \cite{LR} the fibres are non-singular, but
non-singularity is only needed at the intersection with spheres
$\partial B_\epsilon$ and $\partial B_{\epsilon'}$.  In our setting
both fibres are singular Milnor fibres of functions with isolated
singularity on stratified hypersurfaces and in such a case, L\^{e}'s
result \cite{Le} tells that they are, homotopically, wedges of spheres
of dimension $n-1$. Now by \cite[Theorem 3.1, Cor. 3.5]{ST}, the
number of spheres is equal to $\lambda_{p,c}(0)$ and
$\lambda_{p,c}(s)$ respectively. Since $\lambda_{p,c}(s)$ is
independent of $s$ by Proposition \ref{c:dep}, these two numbers
coincide.
 
This shows that one may apply the h-cobordism theorem and conclude
that there exists a $\mathrm{C}^\infty$ function without critical points
on the manifold $(B_\epsilon\setminus \mathring B_{\epsilon'})\cap
g_s^{-1}(c)$, having as levels $\partial B_\epsilon\cap g_s^{-1}(c)$
and $\partial B_{\epsilon'}\cap g_s^{-1}(c)$.  This function can be
extended, with the same property, first on a thin tube
$(B_\epsilon\setminus \mathring B_{\epsilon'})\cap g_s^{-1}(\Delta)$,
where $\Delta$ is a small enough disk centred at $c$, then further
gluing to the distance function on $B_{\epsilon'}\cap
g_s^{-1}(\Delta)$. This extension plays now the role of the distance
function in the construction of the contracting vector field on $T_s$.
Finally, this vector field is used to extend the homeomorphism $h$
from the boundary to the interior of $T_s$, by a similar construction
as the one used in proving condition (\ref{it:l1}).

The conditions (\ref{it:l1}) and (\ref{it:l2}) are now proved and therefore Lemma
\ref{lem:3tim} can be applied. Its conclusion is just Theorem
\ref{th:muinf}.
\end{proof}

\section{Proof of the main theorem}
\label{sec:proof}
  
We first prove Theorem \ref{th:mucst} in case the coefficients of the
family $P$ are polynomials in the variable $s$. The general case of
continuous coefficients will follow by a constructibility argument.

\subsection{Transversality in the neighbourhood of infinity}
\label{ss:trans}

Let $R_1 >0$ such that for all $R\ge R_1$ and all $c\in \Binf(0)$ the
intersection $f_0^{-1}(c) \cap S_R$ is transversal.  We choose $0 <
\eta \ll 1$ such that for all $c\in \Binf(0)$ and all $t\in D_\eta(c)$
the intersection $f_0^{-1}(t) \cap S_{R_1}$ is transversal.  We set
$$K(0) = D \setminus \bigcup_{c\in\Binf(0)} \mathring D_{\eta}(c)$$
for a sufficiently large disk $D$ of $\Cc$.  There exists $R_2 \ge
R_1$ such that for all $t\in K(0)$ and all $R\ge R_2$ the intersection
$f_0^{-1}(t)\cap S_R$ is transversal (see \cite[Prop. 2.11, Cor.
2.12]{Ti-reg} for a more general result, or see \cite[Lemma 5]{Bo}).

By Proposition \ref{c:dep}, $\B(s)$ is an algebraic braid so we may
assume that for a large enough $D$, $\B(s) \subset \mathring D$ for
all $s\in [0,u]$.  Moreover there exists a diffeomorphism $\chi :
\Cc\times [0,u] \to \Cc\times [0,u]$ with $\chi(x,s) = (\chi_s(x),s)$
and such that $\chi_0=\id$, that $\chi_s(\B(s))=\B(0)$ and that
$\chi_s$ is the identity on $\Cc \setminus \mathring D$, for all $s\in
[0,u]$.  We set $K(s) = \chi_s^{-1}(K(0))$.

We may choose $u$ sufficiently small such that for all $s \in [0,u]$,
for all $c\in \Binf(0)$ and all $t\in \chi_s^{-1}(D_\eta(c))$ the
intersection $f_s^{-1}(t) \cap S_{R_1}$ is transversal.  We may also
suppose that for all $s \in [0,u]$, for all $t\in K(s)$ the
intersection $f_s^{-1}(t)\cap S_{R_2}$ is transversal.  Notice that
the intersection $f_u^{-1}(t)\cap S_{R}$ may not be transversal for
all $R \ge R_2$ and $t\in K(s)$.

\subsection{Affine part}
We denote
$$
B'(s) = \big( f_s^{-1}(D) \cap B_{R_1}\big) \cup \big(
f_s^{-1}(K(s)) \cap B_{R_2}\big), \quad s\in [0,u].$$

By using Timourian's theorem at the affine singularities and by gluing
the pieces with vector fields as done in \cite[Lemma 15]{Bo}, we get
the following trivialisation:

$$
\xymatrix{ {} B' \ar[r]^-{\Omega^\mathit{aff}} \ar[d]_-{F}
  &  B'(0) \times [0,u] \ar[d]^-{f_0 \times \id} \\
  D \times [0,u] \ar[r]_-\chi & D \times [0,u], }
$$
where $B'= \bigcup_{s\in [0,u]}B'(s)\times \{s\}$ and $F(x,s) =
(f_s(x),s)$.

\subsection{At infinity, around an atypical value}
It remains to deal with the part at infinity $f_s^{-1}(D) \setminus
\mathring B'(s)$ according to the decomposition of $D$ as the union of
$K(s)$ and of the disks around each $c\in \Binf(s)$.  For each
singular point $(p,c(0))$ at infinity we have a Milnor tube $T_{p,0}$
defined by a Milnor ball of radius $\epsilon(p,c(0))$ and a disk of
radius $\eta$, small enough in order to be a common value for all such
points.

Let $g_s$ be the restriction to $\Xx(s)$ of the compactification
$\tau$ of the $f_s$, let $G : \Xx \longrightarrow \Cc \times [0,u]$ be
defined by $G(x,s)=(g_s(x),s)$ and let $C'(s) =
g_s^{-1}(\chi_s^{-1}(D_\eta(c(0)))) \setminus (\mathring B_{R_1} \cup
\bigcup_{p(s)} \mathring T_{p(s)})$.  Now $g_s$ is transversal to the
following manifolds: to $T_{p(s)}\cap \partial B_\epsilon$, for all
$s\in [0,u]$, by the definition of a Milnor tube, and to $S_{R_1}\cap
C'(s)$, by the definition of $R_1$. We shall call the union of these
sub-spaces the {\em boundary of $C'(s)$}, denoted by $\delta C'(s)$.
Let us recall from \S\ref{sec:def} the definition of the Whitney
stratification on $\Xx(s)$ and remark that $C'(s)\cap \Sing \Xx(s) =
\emptyset$.  Therefore $g_s$ is transversal to the stratum $C'(s)\cap
\Xx^\infty(s)$.

Let $C' = \bigcup_{s\in[0,u]} C'(s)\times \{ s \}$ and remark that $C'
\cap \Sing \Xx = \emptyset$ and that the stratification $\{C'\setminus
\Xx^\infty,\Xx^\infty \}$ is Whitney.  Then by our assumptions and for
small enough $u$, the function $G$ has maximal rank on
$\bigcup_{s\in[0,u]} \mathring C'(s)\times \{ s \}$, on its boundary
$\delta C' = \bigcup_{s\in[0,u]} \delta C'(s)\times \{ s \}$ and on
$C'\cap \Xx^\infty(s)$.  By Thom-Mather's first Isotopy Theorem, $G$ is
a trivial fibration on $C'\setminus \Xx^\infty$. More precisely, one
may construct a vector field on $C'$ which lifts the vector field
$(\frac{\partial \chi_s}{\partial s},1)$ of $D\times [0,u]$ and which
is tangent to the boundary $\delta C'$ and to $C'\cap \Xx^\infty$.  We
may in addition impose the condition that it is tangent to the
sub-variety $g_s^{-1}(\partial \chi_s^{-1}(D_\eta(c(0)))) \cap
S_{R_2}$. We finally get a trivialisation of $C'$, respecting fibres
and compatible with $\chi$.

\subsection{Gluing trivialisations by vector fields}
Since this vector field is constructed such that to coincide at the
common boundaries with the vector field defined on each tube $T$ in
the proof of Theorem \ref{th:muinf}, and with the vector field on $B'$
as defined above, this enables one to glue all the resulting
trivialisations over $[0,u]$. Namely, for
$$B''(s) := \big( f_s^{-1}(D) \cap B_{R_2} \big) \cup \big(
f_s^{-1}(D\setminus \mathring K(s)) \big) \ \mbox{and}\ B'':=
\bigcup_{s\in [0,u]}B''(s)\times \{s\}$$
we get a trivialisation:
$$
\xymatrix{ {} B'' \ar[r]^-{\Omega} \ar[d]_-{F}
  &  B''(0) \times [0,u] \ar[d]^-{f_0 \times \id} \\
  D \times [0,u] \ar[r]_-\chi & D \times [0,u].  }
$$
This diagram proves the topological equivalence of the maps $f_0 :
B''(0) \longrightarrow D$ and $f_u : B''(u) \longrightarrow D$.

\subsection{Extending topological equivalences}
By the transversality of $f_0^{-1}(K(0))$ to the sphere $S_R$, for all
$R\ge R_2$, it follows that the map $f_0 : B''(0) \longrightarrow
\mathring D$ is topologically equivalent to $f_0 : f_0^{-1}(\mathring
D) \longrightarrow \mathring D$, which in turn is topologically
equivalent to $f_0 : \Cc^n \longrightarrow \Cc$.

We take back the argument for $f_0$ in \S \ref{ss:trans} and apply it
to $f_u$: there exists $R_3 \ge R_2$ such that $f_u^{-1}(t)$
intersects transversally $S_R$, for all $t\in K(u)$ and all $R \ge
R_3$.  Now, with arguments similar to the ones used in the proof of
the classical L\^e-Ramanujam theorem (see e.g. \cite[Theorem
5.2]{Ti-reg} or \cite[Lemma 8]{Bo} for details), we show that our
hypothesis of the constancy of $\mu+\lambda$ allows the application of
the h-cobordism theorem on $B'''(u)\setminus \mathring B''(u)$, where
$B'''(u) = \big( f_s^{-1}(D) \cap B_{R_3} \big) \cup \big(
f_s^{-1}(D\setminus \mathring K(s)) \big)$. Consequently, we get a
topological equivalence between $f_u : B''(u) \longrightarrow D$ and
$f_u : B'''(u) \longrightarrow D$.  Finally $f_u : B'''(u)
\longrightarrow \mathring D$ is topologically equivalent to $f_u :
f_u^{-1}(\mathring D) \longrightarrow \mathring D$ by the
transversality evoked above, and this is in turn topologically
equivalent to $f_u : \Cc^n \longrightarrow \Cc$.

\subsection{Continuity of the coefficients}
\label{ss:cont}
So far we have proved Theorem \ref{th:mucst} under the hypothesis that
the coefficients of the family $P$ are polynomials in the parameter
$s$.  We show in the last part of the proof how to recover the case of
continuous coefficients. The following argument was suggested to the
first named author by Frank Loray.  Let $\PP_d(\mu+\lambda, \#\B)$ be
the set of polynomials of degree $d$, with isolated singularities in
the affine space and at infinity, with fixed number of vanishing
cycles $\mu+\lambda$ and with a fixed number of atypical values $\#
\B$.  Proposition \ref{l:construc} tells that $\PP_d(\mu+\lambda,
\#\B)$ is a constructible set.  Since $f_0$ and $f_1$ are in the same
connected component of $\PP_d(\mu+\lambda, \#\B)$, we may connect
$f_0$ to $f_1$ by a family $g_s$ with $g_0=f_0$ and $g_1 = f_1$ such
that the coefficients of $g_s$ are piecewise polynomial functions in
the variable $s$.  Using the proof done before for each polynomial
piece, we finally get that $f_0$ and $f_1$ are topologically
equivalent.  This completes the proof of Theorem \ref{th:mucst}.


\end{document}